\documentclass{article}
\usepackage[utf8]{inputenc}
\usepackage{amsmath}
\usepackage{amsthm}
\usepackage{amsfonts}
\usepackage{tikz}
\usepackage{tabularx}
\usepackage{graphicx}
\usepackage{caption}
\usepackage{pgfplots}
\usepackage{hyperref}
\usepackage{subfig}
\usepackage{multirow}
\usepackage{hhline}
\usepackage[backend=bibtex,giveninits=true]{biblatex}
\usepackage[a4paper, margin=2.5cm]{geometry}
\usepackage[shortlabels]{enumitem}

\title{Global optimality in model predictive control via hidden invariant convexity}
\author{Jorn H. Baayen\footnote{KISTERS Group, Business Unit Water, Utrecht, The Netherlands. \newline E-mail: \url{jorn.baayen@kisters-bv.nl}.} \, and Krzysztof Postek\footnote{Delft University of Technology, Delft, The Netherlands.\newline E-mail: \url{k.s.postek@tudelft.nl}}}
\newtheorem{thm}{Theorem}[section]

\newtheorem{cor}{Corollary}
\newtheorem{defn}{Definition}[section]

\newtheorem{exmp}{Example}[section]

\DeclareMathOperator{\interior}{int}

\bibliography{main}

\begin{document}
\maketitle

\begin{abstract}
Non-convex optimal control problems occurring in, \emph{e.g.}, water or power systems, typically involve a large number of variables related through nonlinear equality constraints. The ideal goal is to find a globally optimal solution, and numerical experience indicates that algorithms aiming for Karush-Kuhn-Tucker points often find (near-)optimal solutions. In our paper, we provide a theoretical underpinning for this phenomenon, showing that on a broad class of problems the objective can be shown to be an \emph{invariantly convex} function (\emph{invex} function, see \cite{craven1981invex,hanson1981sufficiency}) of the control decision variables when state variables are eliminated using implicit function theory. In this way, near-global optimality can be demonstrated, where the exact nature of the global optimality guarantee depends on the position of the solution within the feasible set. In a numerical example, we show how high-quality solutions are obtained with local search for a river control problem where invexity holds.
\end{abstract}

\section{Introduction}
%\red{WHY IMPORTANT}
Model predictive control (MPC) is an ubiquitous technique for optimal control of systems driven by nonlinear PDEs, such as water, gas, and power systems  \cite[e.g.]{ackermann2000real,burgschweiger2009optimization,hante2017challenges}. A typical objective is to steer the system into tracking target levels for the state variables, \emph{e.g.}, stabilizing water level around the desired level \cite{Garcia1989}. %The quantity to minimize is then a norm of the states' deviation from the target. %, plus a control-related cost.

%\red{WANTED BUT DIFFICULT}
One of the key requirements for a solution is \emph{global optimality}. Although locally optimal solutions are often acceptable, globally optimal ones typically yield substantially better objective values and are a more solid base for major commitments. While approaches aiming for global optimality in general nonlinear optimization have been proposed \cite[e.g.]{sahinidis1996baron,belotti2009branching,misener2014antigone,ghaddar2017polynomial}, it remains a challenge in case of tight computation time limits or large problem sizes. Therefore, common work-arounds are to use linearizations \cite[e.g.]{Eschenbach2001,amann2016online,falk2016model} or convex restrictions or relaxations  \cite[e.g.]{madani2014convex,horvath2019convex,lee2019convex} which provide tractability yet at the cost of model accuracy, or to resort to genetic algorithms  \cite[e.g.]{van2004operational,nicklow2010state,vermuyten2018combining}.

It is most desirable, however, to obtain a globally optimal solution to the ``most exact'' nonlinear model without resorting to computationally expensive techniques. Ample numerical evidence exists that local solutions to nonlinear MPC problems are often of high quality, hardly distinguishable from true global optima \cite[e.g.]{burgschweiger2009optimization,ghaddar2017polynomial,Baayen2019-3}. 

In this paper, we provide a theoretical underpinning for this phenomenon by showing that, in the problems we study, the objective function composed with the dynamics is \emph{invex} in the original sense of Hanson \cite{hanson1981sufficiency} and Craven \cite{craven1981invex}. Invexity is a generalization of convexity that certifies, roughly speaking, that on an open set a stationary point of an invex function is its global minimum. Despite the large body of theoretical work, we are not aware of research leveraging invexity in large-scale applied (engineering) context.

We use invexity to prove near-global optimality guarantees for KKT points of optimization problems built for MPC, including full global optimality for KKT points in the interior of the feasible set. Because the invexity of the involved functions is not readily seen and follows from elimination of the problem's constraints, we refer to it as \emph{hidden invexity}.

%\red{IMPLICATION OF OUR RESULTS}
The research contributions are as follows.
\begin{enumerate}[(a)]
\item From the optimal control angle, we show that for a large class of numerical optimal control problems, invex formulations exist that preserve the exact nonlinear dynamics. This allows to tractably determine high-quality solutions to large-scale non-convex MPC problems. 
\item From the mathematical optimization/operations research angle, we show that on difficult problems, certifying the hidden invexity and using standard local search algorithms is a viable alternative to the use of tools designed for general non-convex optimization problems.
\item From the nonlinear analysis angle, we show that invexity of the objective of a constrained optimization problem yields near-global optimality guarantees for KKT points with active inequality constraints, and global optimality for KKT points in the interior of the domain.
\end{enumerate}

The remainder of this paper is structured as follows. Section~\ref{sec:regular.problems} introduces the notion of regular MPC problems for which we establish our result. In Section~\ref{sec:largep} we prove the main result of the paper. Section~\ref{sec:experiment} presents a numerical study for a single river reach modelled using a nonlinear PDE.

\section{Regular MPC problems} \label{sec:regular.problems}
\begin{figure}[h]
\centering
\begin{tikzpicture}

\foreach \m/\l [count=\y] in {1,2}
  \node [circle, fill, minimum size=0.2cm, inner sep=0pt] (input-\m) at (0,1.5-\y) {};

\foreach \m [count=\y] in {1,2,3}
  \node [circle, fill, minimum size=0.2cm, inner sep=0pt] (hidden-\m) at (2,2.5-\y*1.25) {};

\foreach \m [count=\y] in {1,2}
  \node [circle, fill, minimum size=0.2cm, inner sep=0pt] (output-\m) at (4,1.5-\y) {};
  
\node [circle, fill, minimum size=0.2cm, inner sep=0pt] (objective) at (6,0.0) {};

\foreach \l [count=\i] in {1,2}
  \draw [<-] (input-\i) -- ++(-1,0)
    node [above, midway] {$u_\l$};

\foreach \l [count=\i] in {1,2,3}
  \node [above] at (hidden-\i.north) {$x_\l$};

\foreach \l [count=\i] in {1,2}
  \draw [->] (output-\i) -- (objective);

\foreach \l [count=\i] in {1,2}
  \node [above] at (output-\i.north) {$y_\l$};   
  
\node [above] at (objective.north) {$f$};

\foreach \i in {1,...,2}
  \foreach \j in {1,...,3}
    \draw [->] (input-\i) -- (hidden-\j);

\foreach \i in {1,...,3}
  \foreach \j in {1,...,2}
    \draw [->] (hidden-\i) -- (output-\j);

\end{tikzpicture}
\captionsetup{font=scriptsize}
\caption{The relationship between the controls $u$, the implicitly defined states $x$, and the output variables $y$, for $n=2$ and $m=3$.}
\label{fig:schema}
\end{figure}
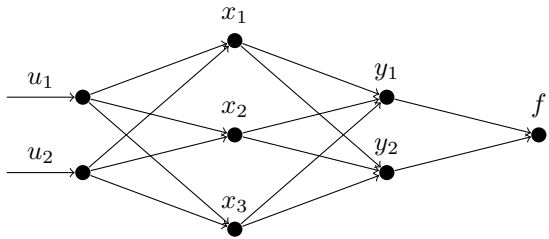

% \subsection{Definitions}
In this section we describe the class of optimization problems for which we demonstrate invexity. Consider the optimization problem
\begin{align}\label{eq:general-prob}
\min_{x,u}\, (f \circ g)(x) & \quad \text{subject to} \tag{$\mathcal{P}$} \\
c(x,u) & = 0 \nonumber \\
d(u) & \leq 0, \nonumber
\end{align}
where we refer to the variables $x \in \mathbb{R}^m$ as \emph{states} and the variables $u \in \mathbb{R}^n$ as \emph{controls}, which is because the values of the controls \emph{implicitly} determine the values of the states through the \emph{equality constraints} $c(x,u)=0$. The function $f: \mathbb{R}^n \to \mathbb{R}$ is the \emph{objective} and the function $g: \mathbb{R}^m \to \mathbb{R}^n$ is the \emph{output function} mapping states $x$ to \emph{outputs} $y:=g(x)$.  The relationship between the controls $u$, the implicitly defined states $x$, and the output variables $y$, is illustrated in Figure \ref{fig:schema}.

We denote the set of of admissible controls as $U := \{u \in \mathbb{R}^n: \ d(u) \leq 0 \}$, where $d$ are the \emph{bound inequality constraints}, and denote the set of indices of the coordinates of $d$ with $\mathcal{I}$. 

Our goal will be to show that the objective of problem~\ref{eq:general-prob} is invex as a function of the controls $u$, under certain conditions. However, the equality constraints in \eqref{eq:general-prob} can involve nonlinear functions, making the analysis cumbersome. We shall alleviate this difficulty by eliminating the constraints using implicit function theory and analyzing the problem using \emph{total} gradients with respect to $u$, wherein the derivatives of the state variables $x_i$ with respect to $u$ are expressed explicitly.  This step is used for the analysis, but is not required in practice.

If the Jacobian $\nabla_x c$ is invertible, then the total Jacobian of the states $x$ with respect to the controls $u$ may be expressed using the implicit function theorem as
\begin{equation}\label{eq:total-grad}
D_u x = -\nabla^{-1}_x c \nabla_u c,
\end{equation}
in which the prefix $\nabla_x$ denotes the matrix of partial derivatives with respect to the components of $x$, and $D_x$ the matrix of of total derivatives with respect to $x$, of a given function. 

In order to establish our result, we need some assumptions. In the following, we define the  \emph{regular MPC problems} for which invexity can be demonstrated. After the definition, we discuss each of the conditions, most of which are typical of MPC problems. %\pagebreak

\begin{defn}\label{def:rmpc}
Consider an optimization problem \eqref{eq:general-prob}. Let the functions $f$, $g$, $c$ and $d$ be continuously differentiable. We say that \eqref{eq:general-prob} is a \emph{regular} MPC (rMPC) problem if the following conditions are satisfied:
\begin{enumerate}
\item the set of admissible controls $U$ is bounded using the inequality constraints
\begin{equation*}\label{eq:ctrl-bound}
d(u) \leq 0,
\end{equation*}
where the inequality depends on the controls $u$ only and holds component-wise,
\item no explicit constraints or bounds are imposed on the states $x \in \mathbb{R}^m$,
\end{enumerate}
as well as \emph{linear independence constraint qualifications} (LICQ):
\begin{enumerate}
\setcounter{enumi}{2}
\item the Jacobian matrix of the equality constraints $c$ with respect to the state variables $x$, \emph{i.e.}, $\nabla_x c(x,u)$, is square and full-rank for all $(x,u)$ such that $c(x,u)=0$, $u \in U$,
\item the Jacobian matrix of the equality constraints $c$ with respect to the control variables $u$, \emph{i.e.}, $\nabla_u c(x,u)$, is full-rank for all $(x,u)$ such that $c(x,u)=0$, $u \in U$,
\item the gradient vectors of the active inequality constraints $d_i$ at the point $u$, \emph{i.e.}, $\nabla_u d_i(u)$ for all $i \in \mathcal{I}$ such that $d_i(u)=0$, are linearly independent for all $u \in U$,
\end{enumerate}
the \emph{uniqueness condition}:
\begin{enumerate}
\setcounter{enumi}{5}
\item for every $u \in U$, the constraints $c(x,u)=0$ have a unique solution $x$,
\end{enumerate}
the \emph{output controllability condition}:
\begin{enumerate}
\setcounter{enumi}{6}
\item the output function $g: \mathbb{R}^m \to \mathbb{R}^n$ is such that the square matrix
\begin{equation}\label{eq:invertible}
-\nabla_x g (x) \, \nabla^{-1}_x c(x,u) \, \nabla_{u} c(x,u) 
\end{equation}
is invertible for all $(x,u)$ such that $c(x,u)=0$, $u \in U$,
\end{enumerate}
and the \emph{convexity condition}:
\begin{enumerate}
\setcounter{enumi}{7}
\item the objective function $f: \mathbb{R}^n \to \mathbb{R}$ is convex.
\end{enumerate}
\end{defn}

We now discuss the meaning of each of the respective assumptions. 

Condition 1 states that the set of feasible controls must be described using a finite number of continuously differentiable inequality constraints.  Standard convex domains such as a ball or a box satisfy this condition. 

Condition 2 is sufficient for the constraints $c(x,u)=0$ to have a solution for every $u \in U$.

Condition 3 is a linear independence constraint qualification (LICQ, \cite[e.g.]{nocedal2006numerical}). It is required for the states $x_i$ to be well defined as implicit functions of $u$, and for us to be able to apply implicit function theory.  Condition 3 is typically straightforward to demonstrate if the dynamics of the underlying model are integrable in time. For the dynamics to be uniquely integrable in time, it is required that the number of states be equal to the number of equations, and furthermore that the Jacobian of the equations with respect to the states be non-singular.   

Condition 4 is also a linear independence constraint qualification that essentially states that at all times, the dynamics are sensitive to all controls. 

Condition 5 is the third LICQ and is satisfied automatically for standard domains such as balls and boxes.

Condition 6 states that the functions $u \mapsto x_i$ are uniquely defined on $U$. It is a standard MPC assumption as the required uniqueness property typically follows if the dynamics are (uniquely) integrable in time.

Condition 7 states that different attainable states should map to different outputs.  Since Conditions 3 and 4 imply that the implicit function $u \mapsto x$ is injective, the range of attainable state values is a subset of $\mathbb{R}^m$ with dimension $n$, on which invertible mappings assigning coordinates in  $\mathbb{R}^{n}$ arise naturally.

Another interpretation of Condition 7 is the following.  If we would require the LICQ and uniqueness conditions to hold on $\mathbb{R}^n$ (rather than  on $U$ only), Condition 7 would imply that for every possible value of the output $y$, there exists a control input $u \in \mathbb{R}^n$ that realizes the output $y$.
It is therefore a type of \emph{output-controllability} \cite[e.g.]{ogata1997modern} condition.

Condition 8 is standard and includes objectives such as $p$-norms raised to the $p$th power with $p \geq 2$.

% \subsection{Non-rMPC constraints}
Before looking at a more interesting example of rMPC problems, it is instructive to consider a few irregular problems to show that the rMPC conditions indeed eliminate some of the well-known NP-hard problems.

\begin{exmp}\label{ex:irregular1}
Let $[0,1] \subset U \subset \mathbb{R}$. If a problem contains a binary-restriction constraint $u(1-u)=0$, then Condition 6 is not satisfied.
\end{exmp}

\begin{exmp}\label{ex:irregular2}
Let $[-1,1] \subset U \subset \mathbb{R}$. If a problem contains a sinusoidal constraint $u=\sin x$, then Condition 6 is, in general (\emph{i.e.}, barring additional structure), not satisfied.
\end{exmp}

\begin{exmp}\label{ex:irregular3}
Let $0 \in U \subset \mathbb{R}$. If a problem contains a bilinear constraint of the form $u=x_1 x_2$, then the LICQ Condition 3 is, in general (\emph{i.e.}, barring additional structure), not satisfied.
\end{exmp}

\begin{exmp}\label{ex:irregular4}
Let $(0,0) \in U \subset \mathbb{R}^2$. If a problem contains a bilinear constraint of the form $x=u_1 u_2$, then the LICQ Condition 4 is, in general (\emph{i.e.}, barring additional structure), not satisfied.
\end{exmp}

We now move to introduce examples of systems meeting the rMPC conditions. One such simple example is if the constraints are linear-affine functions satisfying the appropriate rank conditions.
\begin{exmp}
Consider an MPC problem with control vector $u \in \mathbb{R}^n$, state vector $x \in \mathbb{R}^m$, output vector $y \in \mathbb{R}^n$, and trajectory tracking objective
$$
f(y)=\sum_{i=1}^n |y_i - y_i^t|^p
$$
with $p \in [2, \infty)$ and linear-affine output function
$$
y=Cx+c
$$
subject to the bounds
$$
-\infty < u^L_j \leq u_j \leq u_j^U < \infty \quad j \in \mathbb{N}_n
$$
and linear-affine constraints
$$
c(x,u)=Ax+Bu+b
$$
with matrix $A$ square and invertible, matrices $B$ and $C$ full rank and matrix $C$ such that the square matrix $CA^{-1}B$ is invertible. This is a regular problem.
\end{exmp}
Our next example considers trigonometric constraints, which commonly arise in control of systems with axes of rotation such as vehicles, ships, and aircraft \cite[e.g.]{stengel2015flight}.
\begin{exmp}
An MPC problem with control vector $u \in \mathbb{R}^2$, state vector $x \in \mathbb{R}^2$, output vector $y \in \mathbb{R}^2$, and trajectory tracking objective
$$
f(y)=\sum_{i \in \{1,2\}} |y_i - y_i^t|^p
$$
with $p \in [2, \infty)$ and output function
\[
y=x
\]
subject to the bounds
\begin{align*}
0 < u^L_1 & \leq u_1 \leq u_1^U < \infty \\
0 & \leq u_2 < 2\pi
\end{align*}
and constraints
\begin{align*}
x_1&=u_1\cos u_2, \\
x_2&=u_1\sin u_2,
\end{align*}
is regular.
\end{exmp}
Since bilinear constraints are very common, we also show how many of them can satisfy the rMPC assumptions.
\begin{exmp}\label{ex:bilinear}
An MPC problem with control vector $u \in \mathbb{R}^2$, state vectors $x \in \mathbb{R}^2$ and $z \in \mathbb{R}^2$, output vector $y \in \mathbb{R}^2$, and objective
$$
f(y)=\sum_{i \in \{1,2\} }|y_i - y_i^t|^p
$$
with $p \in [2, \infty)$ and output function
\[
y=z
\]
subject to the bounds
$$
0 \leq u^L_j \leq u_j \leq u^U_j < \infty \quad j \in \{1,2\}
$$
and constraints
\begin{align*}
%u_1 & = a x_1 , \\
u_1 & = x_1 z_0, \\
x_1 & = z_0 - z_1, \\
u_2 & = x_2 z_1, \\
x_2 & = z_1 - z_2,
\end{align*}
with the fixed initial condition $z_0 \in \mathbb{R}$, is regular as long as $z_0$ is chosen such that $z_i \neq 0$, $i \in \{1,2\}$, for all feasible $u$. Conditions 1--5 and 7--8 are readily verified. To verify Condition 6, \emph{i.e.}, that for any $u \in U$ the constraints admit a unique solution, note that the constraints may be solved in the displayed order, starting from the fixed value $z_0$.
\end{exmp}
Situations like Example \ref{ex:bilinear} commonly occur when modelling the generation of a hydroelectric turbine in a power station. Instantaneous generation ($u$) is non-negative and bounded, and it is bilinear in flow ($x$) and the water level difference ($z$) across a dam, which is never zero. At the same time, an increase in flow results ($x$) in a decrease of the water level difference ($z$).  Similar reasoning applies to the power consumption of pumps.

Conditions 1--2, 5, and 7--8 may be satisfied by design.   Conditions 3--4 and 6 are also satisfied by appropriate discretizations of certain hyperbolic PDEs, used to model the following examples.
\begin{exmp}\label{eq:hyper1}
A river or canal network modelled using the Saint-Venant equations, with control authority exercised at weirs, dams, gates, and pumps \cite[e.g.]{baayen2019overview}.
\end{exmp}
\begin{exmp}\label{eq:hyper2}
A drinking water distribution network modelled using the Darcy-Weisbach or Hazen-Williams equations, with control authority exercised at valves and pumps \cite[e.g.]{burgschweiger2009optimization}.
\end{exmp}
\begin{exmp}\label{eq:hyper4}
A natural gas and/or hydrogen distribution network modelled using the isothermal Euler equations, with control authority exercised at valves and compressors \cite[e.g.]{hante2017challenges}.
\end{exmp}
Such PDEs have a time dimension along which they may be integrated, starting from a fixed initial condition, analogous to Example \ref{ex:bilinear}. If the discretization is implicit or semi-implicit in time (in the sense of the implicit Euler method \cite[e.g.]{isaacson1966analysis} and Example \ref{ex:bilinear}), then Conditions 3 and 6 are satisfied if, given values for time step $t$, the Jacobian matrices arising when solving for time step $t+1$ are multi-diagonal. Such discretizations exist for Examples \ref{eq:hyper1}--\ref{eq:hyper4}. Proofs of this fact may be found in \cite[e.g.]{casulli1990semi,casulli1998conservative,todini1988gradient,simpson2011jacobian,hante2017challenges}.  

Conditions 4 and 7 also arise naturally in the context of systems driven by hyperbolic PDEs. Since the number of output variables needs to be equal to the number of controls, for every control we can make the output function $g$ select a state corresponding to the discretization node upon which the control variable acts, or corresponding to a spatially adjacent node if the grid is staggered \cite[e.g.]{Stelling2003}. In a river control example, that would mean that for every dam in a river, the corresponding output variable would be the water level directly upstream of the dam. Such an output function $g$ is linear and injective on the set of attainable states \footnote{Let $t$ denote the first time step at which two control strategies differ. Since for hyperbolic PDEs, disturbances travel at finite velocity \cite[e.g.]{evans2010partial}, the resulting change in state at or adjacent to a perturbed control variable at time $t$ cannot be compensated for by a) perturbations to control variables at times $t' > t$, or b) by any of the other control variables acting at another node in the spatial discretization,  acting at time $t$. Therefore different control strategies produce different output.}, whence it satisfies the LICQ Condition 4 and the invertibility Condition 7.  The function of the remaining states $x_i$ is to ensure physically accurate wave propagation, by means of a sufficiently fine spatial discretization, in between of the control nodes.  Experience shows that such a setup produces appropriate control strategies, as will also be illustrated with a concrete example in Section \ref{sec:experiment}.

\section{Hidden invexity}\label{sec:largep}
\subsection{Introduction and the main result}
In this section we present our main result that rMPC problems have \emph{hidden invexity} when reduced to optimization over control variables. We begin by recalling the definition of invexity.
\begin{defn} \label{def:invexity}
A function $f: X \subset \mathbb{R}^n \to \mathbb{R}$ is \emph{invex} on the open set $X$ if there exists a vector function $\eta(x_2,x_1): X \times X \to \mathbb{R}^n$ such that
\begin{equation}\label{eq:invex}
f(x_2)-f(x_1) \geq \eta^T(x_2,x_1) \nabla f(x_1)
\end{equation}
for all $x_1, x_2 \in X$.
\end{defn}
The name invex follows from \emph{invariant convex} \cite{craven1981invex}.  A function is invex if and only if every stationary point is a global minimum. To see the first implication, set $\nabla f=0$ in Equation \eqref{eq:invex}.  A concise proof of the reverse implication may be found in  \cite{ben1986invexity}. 

The definition of invexity is usually stated for functions defined on open sets, wheareas our goal is to optimize over a closed set $U$. There exists an entire family of extensions of the notion of invexity to constrained optimization problems (KT-invexity \cite{martin1985essence}, HC-invexity \cite{hanson1981sufficiency,craven1981invex,martin1985essence}, Type I/Type II invexity \cite{hanson1999invexity}). However, each of them is difficult to apply to real-world problems like ours, due to the need to find a common function $\eta$ for the objective and the constraints. Instead, we shall stay with the standard notion of invexity and, eliminating the equality constraints from the problem, show invexity of the objective function on the interior of the search space in terms of the control variables. In the process, we extend the analysis to problems with inequality-constrained control variables by analyzing the geometry of the boundary of the search space using tangent cones.

In our analysis, we will use the fact that invexity of functions arises naturally in the composition of convex functions with transformations that are full-rank, \emph{i.e.}, that have an invertible Jacobian  \cite[e.g.]{craven1981duality}. We will now show how rMPC problems fit this scheme.

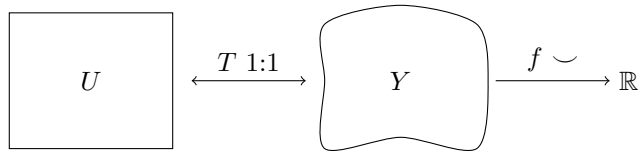
\begin{figure}[h]
\centering
\begin{tikzpicture}

\draw[] (-.15, -.15) rectangle (2.0, 1.65) node[pos=.5] {$U$};

\draw plot [smooth cycle] coordinates {(4.0, -.15)(5.0, 0.0)(6.0, -0.15)(6.15, 0.75)(6.0,1.5)(5.0,1.75)(4.0, 1.5)(4.0, 0.75)} node at (5.0, 0.75) {$Y$};

\node [] (blah) at (8.0, 0.75) {$\mathbb{R}$}; 

\draw [<->] (2.25, 0.75) -- node[above] {$T$ 1:1} (3.75, 0.75);

\draw [->] (6.25, 0.75) -- node[above] {$f$ $\smile$} (blah);

\end{tikzpicture}
\captionsetup{font=scriptsize}
\caption{Convex objective function $f$ composed with invertible transformation $T$.}
\label{fig:transform}
\end{figure}

In our problem we have a convex objective function $f$. Next, we need a full-rank transformation. We construct it by noticing that Conditions 3 and 6 in Definition \ref{def:rmpc} enable us to use the implicit function theorem to express the state variables $x$ as a function $u \mapsto x$. Problem \eqref{eq:general-prob} can therefore be rewritten as:
\begin{align}\label{eq:eliminated}
\min_{u} (f \circ g \circ x)(u) & \quad \text{subject to} \tag{$\mathcal{P}^U$} \\
d(u) & \leq 0. \nonumber
\end{align}
In \eqref{eq:eliminated}, the composition $T: U \to Y:=T[U]$, $T(u):=(g \circ x)(u)$, will be playing the role of the invertible transformation, and the composition $f \circ T$ will be shown to be invex.  This setup is illustrated in Figures  \ref{fig:schema} and \ref{fig:transform}.  

The key feature to deal with in our analysis is the fact that in general, the set $Y$ is non-convex. The strength of our main result for a particular KKT point with $y^*$ will depend on the place where point $y^*$ is in the set $Y$ - in the interior or on the boundary. To make this distinction rigorous and state our main result, we first recall the definition of the tangent cone \cite[e.g.]{geiger2013theorie,nocedal2006numerical}.

\begin{defn}\label{def:tang-cone}
Let $Y \subset \mathbb{R}^n$ be a non-empty set. A vector $d \in \mathbb{R}^n$ is \emph{tangent} to $Y$ at $y \in Y$, if there exist sequences $\{y^k\} \subset Y$, $\{t_k\} \subset \mathbb{R}^+$ such that
$$
y^k \to y, \quad t_k \to 0, \quad  \frac{y^k - y}{t_k} \to d.
$$
The set of all tangent vectors at $y \in Y$ is the \emph{tangent cone} of $Y$ at $y$, denoted $\mathcal{T}_Y(y)$.
\end{defn}

We now state our main result the proof of which is relegated to Section \ref{sec:proof}.

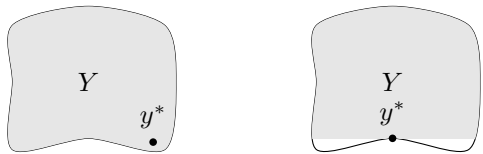
\begin{figure}[h]
\centering
\begin{tikzpicture}

\foreach \x in {0,4}
  \draw plot [smooth cycle] coordinates {(\x, -.15)(\x+1, 0.0)(\x+2, -0.15)(\x+2.15, 0.75)(\x+2,1.5)(\x+1,1.75)(\x, 1.5)(\x, 0.75)} node at (\x + 1, 0.75) {};
  
\foreach \x in {0}
  \fill [gray!20] plot [smooth cycle] coordinates {(\x, -.15)(\x+1, 0.0)(\x+2, -0.15)(\x+2.15, 0.75)(\x+2,1.5)(\x+1,1.75)(\x, 1.5)(\x, 0.75)} node at (\x + 1, 0.75) {};
\begin{scope}
\clip(3,0) rectangle (7,2.0);
\foreach \x in {4}
  \fill [gray!20] plot [smooth cycle] coordinates {(\x, -.15)(\x+1, 0.0)(\x+2, -0.15)(\x+2.15, 0.75)(\x+2,1.5)(\x+1,1.75)(\x, 1.5)(\x, 0.75)} node at (\x + 1, 0.75) {};
\end{scope}

\node [circle, fill, minimum size=0.1cm, inner sep=0pt] (blah) at (1.85, -.05) {}; 
\node [above] at (blah.north) {$y^*$}; 
% Name of the set Y
\node (blah2) at (1.0, 0.75) {}; 
\node at (blah2) {$Y$}; 

\node [circle, fill, minimum size=0.1cm, inner sep=0pt] (blah) at (5.0, 0.0) {}; 
\node [above] at (blah.north) {$y^*$}; 
% Name of the set Y
\node  (blah3) at (5.0, 0.75) {}; 
\node at (blah3) {$Y$}; 

\end{tikzpicture}
\captionsetup{font=scriptsize}
\caption{The highlighted areas illustrate sets $T[V(u^*)]$ within which a solution $y^*=T(u^*)$ is provably globally optimal, for an interior (left) and a boundary solution (right). The sets are shown in the output space $Y=T[U]$ to highlight the role of the tangent cones. Note that in the left panel the point $y^*$ is in the interior of $Y$, hence a global optimum in $Y$, despite the fact that there are some points in $Y$ that are not ``visible'' from the tangent cone at $y^*$. This is thanks to invexity.}
\label{fig:validity}
\end{figure}

\begin{thm}\label{thm:invex}
Consider an rMPC problem \eqref{eq:eliminated}.
Let $(u^*,\lambda^*)$ be a KKT point of this problem. 
Then $u^*$ is a global minimum of $f \circ T$ on the set
\begin{equation}\label{eq:vset}
V(u^*):=\{u \in U : T(u) - T(u^*) \in \mathcal{T}_Y(T(u^*))\},
\end{equation}
where $\mathcal{T}_Y(T(u^*))$ denotes the tangent cone of $Y$ at $T(u^*)$.
\end{thm}
It states that, in rough terms, a KKT point is a global optimum with respect to the interior of the domain and all inactive boundary segments (invexity), minus any points ``hidden from view'' due to local non-convexity of the active boundary segments. The geometric meaning is illustrated in Figure \ref{fig:validity} and explored in further detail with the corollaries below and subsequent discussion.

Before proceeding to the corollaries, we note that the reverse statement of Theorem \ref{thm:invex}, \emph{i.e.}, that every minimum is a KKT point, follows   from the LICQ Conditions 3--5 in Definition \ref{def:rmpc}. The LICQ conditions form the regularity condition required for every minimum to be a KKT point    \cite[e.g.]{nocedal2006numerical}.

The first corollary is a direct consequence of the fact that for an interior point, $\mathcal{T}_Y(T(u^*))=\mathbb{R}^n$.
\begin{cor}
Consider an rMPC problem \eqref{eq:eliminated}. Consider a KKT point $(u^*, \lambda^*)$ such that $u^* \in \interior U$. Then $u^*$ is a global minimum of $f \circ T$ on $U$.
\end{cor}

In other words, $f \circ T$ is invex on the interior of $U$.  

The second corollary follows from the fact that for a convex set $X$, membership $x, y \in X$ implies that $y - x \in \mathcal{T}_X(x)$. % and vice versa.

\begin{cor}\label{cor:convex}
Consider an rMPC problem \eqref{eq:eliminated}. Consider a KKT point $(u^*, \lambda^*)$ and a set $W \subset U$, $u^* \in W$, such that $T[W]$ is convex.  Then $u^*$ is a global minimum of $f \circ T$ on the set $W$.
\end{cor}

We will now explain the meaning of these results. For this, it is instructive to first recall the reference situation: general nonlinear programming. A KKT point of a nonlinear optimization problem need not be a local minimum; it may also be a local maximum, or a saddle point.  Furthermore, in case that a KKT point is a local minimum, it is only guaranteed to be minimal within an arbitrarily small neighbourhood of itself. From a numerical point of view, generic nonlinear optimization problems are hard: local search methods may converge to KKT points that are local maxima or saddle points.

For an rMPC problem, Theorem \ref{thm:invex} provides a stronger characterization of KKT points. First of all, it states that a KKT point is a local minimum, and therefore not a maximum or a saddle point, within $V(u^*)$. This is important from a numerical point of view.  Secondly, it states that a local minimum $u^*$ is a global minimum within the set $V(u^*)$. If $u^*$ is an interior point, we have that $V(u^*)=U$, whence it is a global optimum.  If $u^*$ lies on the boundary, its objective value is no greater than the objective values for all points that corresponding to points that lie on rays emanating from $y^* = T(u^*)$ in the directions of the tangent vectors in set $Y$, a situation we refer to as \emph{near-globally optimal}.  Both cases are illustrated in Figure \ref{fig:validity}.

The term \emph{hidden invexity} of an rMPC problem arises exactly from the need to analyze the behaviour of the composition of a convex objective function with an invertible transformation. If, furthermore, the set $Y$ is convex, then by Corollary \ref{cor:convex} the rMPC problem may be said to have \emph{hidden convexity} \cite[e.g.]{ben1996hidden,li2005hidden}.

\subsection{Proof of Theorem~\ref{thm:invex}}\label{sec:proof}

Consider the transformation $T = g \circ x$. By Condition 7 of Definition \ref{def:rmpc}, $D_u T$ is invertible whence, by the inverse function theorem, $T$ itself is invertible.  The transformation and its use within the optimization problem is illustrated in Figure \ref{fig:transform}.

We will first show that a point $(u^*,\lambda^*)$ is a KKT point of the optimization problem 
\begin{align}\label{eq:optU}
\min_{u} \, (f \circ T)(u) & \quad  \text{subject to} \tag{$\mathcal{P}^U$} \\ 
d(u) & \leq 0, \nonumber
\end{align}
if and only if $(T(u^*),\lambda^*)$ is a KKT point of the optimization problem
\begin{align}\label{eq:optY}
\min_{y} f(y) & \quad  \text{subject to} \tag{$\mathcal{P}^Y$} \\ 
(d \circ T^{-1})(y) & \leq 0, \nonumber
\end{align}
Afterwards, we will analyze the global optimality structure of the KKT points.

Let
$$\mathcal{L}^U(u,\lambda):=(f \circ T)(u)+\lambda^T d(u)$$ denote the Lagrangian of problem  \eqref{eq:optU}, and let 
$$\mathcal{L}^Y(y,\lambda):=f(y)+\lambda^T (d \circ T^{-1})(y)$$ denote the Lagrangian of problem  \eqref{eq:optY}. We will use the standard definition of KKT points following \cite{nocedal2006numerical}.
KKT points of \eqref{eq:optU} are stationary points of the Lagrangian $\mathcal{L}^Y$  and therefore satisfy
\begin{align*}
0 & = D_y \mathcal{L}^Y \\
& = \nabla_y f + \lambda^T \nabla_u d \, D^{-1}_u T \\
& = \nabla_y f \left[D_u T \, D^{-1}_u T\right] + \lambda^T \nabla_u d \, D^{-1}_u T \\
& = \left[\nabla_y f \, D_u T + \lambda^T \nabla_u d\right] D^{-1}_u T \\
& = D_u \mathcal{L}^U \, D^{-1}_u T.
\end{align*}
Since $D_u T$ is invertible, a point $(u^*,\lambda^*)$ is a stationary point of $\mathcal{L}^U$ if and only if $(T(u^*),\lambda^*)$ is a stationary point of $\mathcal{L}^Y$.  Similar reasoning applies to the primal and dual feasibility conditions ($d(u^*) \leq 0$ and $\lambda^* \geq 0$) as well as to the complementarity condition ($\lambda^*_i d_i(u^*)=0$ whenever $d_i(u^*)=0$). This completes the first part of the proof.

We will now analyze the KKT points.  For this, rather than using the definition of invexity directly, we will use some machinery related to tangent cones.  In this way, we will also be able to reason about points on the boundary of the search space; recall that invexity is defined on open sets, \emph{i.e.}, sets without their boundary (cf. Definition \ref{def:invexity}). 

We first recall a few additional definitions.  Relevant references are \cite{aubin1984applied,geiger2013theorie,nocedal2006numerical}.

\begin{defn}
The set
$$
\mathcal{A}(y^*):=\{i \in \mathcal{I}: (d_i \circ T^{-1})(y^*)=0 \}
$$
is the \emph{active set} for the problem \eqref{eq:optY} at the point $y^* \in Y$.
\end{defn}

\begin{defn}
The set
$$
\mathcal{F}(y^*):=\{t \in \mathbb{R}^n : t^T D_y (d_i \circ T^{-1})(y^*) \leq 0 \quad \forall i \in \mathcal{A}(y^*)\}, 
$$
is the \emph{set of linearized feasible directions} for the problem \eqref{eq:optY} at the point $y^* \in Y$.
\end{defn}

\begin{defn}
The cone
$$
K^{\circ}:=\{y \in \mathbb{R}^n : y^T x \leq 0 \quad \forall x \in K \}
$$
is the \emph{polar cone} of the cone $K$.
\end{defn}

Let $(u^*,\lambda^*)$ be a KKT point of \eqref{eq:optU}.  Our aim is to show that the point $u^*$ is global minimum of $f$ on the set $V(u^*)$ as defined in Equation \eqref{eq:vset}. For this, it is convenient to reason about $y^*=T(u^*)$ and problem \eqref{eq:optY}. By virtue of Condition 5 of Definition \ref{def:rmpc} LICQ holds for the constraint function $d \circ T^{-1}$, whence $\mathcal{F}(y^*)=\mathcal{T}_Y(y^*)$. See \cite[e.g.]{nocedal2006numerical} for proof of this fact. 

Since, by the first part of this proof, $(y^*,\lambda^*)$ is also KKT point, we have
\begin{equation}\label{eq:kkt-eq}
-\nabla_y f(y^*)=\lambda^{*T} D_y (d \circ T^{-1})(y^*).
\end{equation}
Following the definition of the set of linearized feasible directions $\mathcal{F}(y^*)$, for all $t \in \mathcal{F}(y^*)=\mathcal{T}_Y(y^*)$ we have that $t^T D_y (d_i \circ T^{-1})(y^*) \leq 0$ for all $i \in \mathcal{A}(y^*)$. Because of this and the facts that $\lambda^*_i \geq 0$ for all $i \in \mathcal{A}(y^*)$ and $\lambda^*_i=0$ for all $i \in \mathcal{I} \setminus \mathcal{A}(y^*)$, it follows from Equation \eqref{eq:kkt-eq} that $-t^T \nabla_y f(y^*) \leq 0$ for all $t \in \mathcal{F}(y^*)=\mathcal{T}_Y(y^*)$. Therefore $-\nabla_y f(y^*) \in (\mathcal{T}_Y(y^*))^{\circ}$, the polar cone of the tangent cone.  

Since $T[V(u^*)] \subset Y$, it follows directly from Definition \ref{def:tang-cone} that $\mathcal{T}_{T[V(u^*)]}(y^*) \subset  \mathcal{T}_Y(y^*)$. The inclusion reverses when taking polar cones, so that $$-\nabla_y f(y^*) \in (\mathcal{T}_Y(y^*))^{\circ} \subset (\mathcal{T}_{T[V(u^*)]}(y^*))^{\circ}.$$In other words, for every tangent vector $t \in \mathcal{T}_{[V(u^*)]}(y^*)$, we have $t^T \nabla_y f(y^*) \geq 0$.   By convexity of $f$, for every $y \in T[V(u^*)]$, $$f(y)-f(y^*) \geq (y-y^*)^T \nabla_y f(y^*) \geq 0.$$ The second inequality follows from the fact that $y-y^* \in \mathcal{T}_{Y}(y^*)$ by construction of the set $V(u^*)$.  We conclude that $y^*$ is a global minimum of $f$ on $T[V(u^*)]$, whence $u^*$ is a global minimum of $f \circ T$ on $V(u^*)$.  \qed

\section{Numerical experiment}\label{sec:experiment}
In this section, we describe a numerical experiment revolving around an MPC problem for the one-dimensional \emph{shallow water equations}. The experiment illustrates how local search indeed leads to near-globally optimal solutions in the sense of this paper.  This furthermore highlights the practical relevance of our result, since, traditionally, one would only be able to claim (near-)global optimality when using computationally expensive general-purpose solvers for non-convex problems such as, \emph{e.g.}, \textsc{Couenne} \cite{belotti2009branching}, or, alternatively, seeding a local search with a large number of different starting points in order to obtain increased confidence in the quality of the solution.

The one-dimensional shallow water equations are also known as the  \emph{Saint-Venant equations}, and form a nonlinear hyperbolic PDE  \cite[e.g.]{vreugdenhil2013numerical}.  
The Saint-Venant equations describe levels and flows in rivers and canals. They are given by the momentum equation
\begin{equation*}
\frac{\partial Q}{\partial t} + \frac{\partial}{\partial x}\frac{Q^2}{A} + gA\frac{\partial H}{\partial x} + g\frac{Q|Q|}{A R C^2}=0,
\end{equation*}
with longitudinal coordinate $x$, time $t$, discharge $Q$, water level $H$, cross section $A$, hydraulic radius $R:=A/P$, wetted perimeter $P$, Ch\'{e}zy friction coefficient $C$, gravitational constant $g$, and by the mass balance (or continuity) equation
\begin{equation*}
\frac{\partial Q}{\partial x} + \frac{\partial A}{\partial t} = 0.
\end{equation*}
In a setting with bidirectional flow, the $|Q|$ factor in the momentum equation may be approximated by a smooth function \cite[e.g.]{burgschweiger2009optimization}.  In this section, however, we will only consider unidirectional flow with $Q > 0$, so that $|Q|=Q$.

For our experiment, we consider a discretization that is semi-implicit in time, following \cite{casulli1998conservative}, and staggered in space. The semi-implicit discretization ensures that Conditions 3--4 and 6 in Definition \ref{def:rmpc} are met. The proof of this statement is given, in a different wording, in \cite[e.g.]{casulli1990semi,casulli1998conservative}. We use the experimental setting from the draft  \cite{baayen2020continuation}, from which the following description and Table \ref{table:example-specs} are  adapted.  

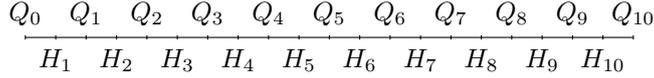
\begin{figure}[!t]
\centering
\begin{tikzpicture}[scale=0.8]
\draw (0,0) -- (10.0,0);
\foreach \x in {0,1,...,10}
    \draw (\x cm,1 pt) -- (\x cm,-1 pt) node[anchor=south] { \begin{tabular}{c} $Q_{\x}$ \end{tabular} };
\foreach \x in {1,2,...,10}
    \draw (\x cm - 0.5 cm,1 pt) -- (\x cm - 0.5 cm,-1 pt) node[anchor=north] {$H_{\x}$};
\end{tikzpicture}
\captionsetup{font=scriptsize}
\caption{Staggered grid for the example problem.}
\label{fig:example-grid}
\end{figure}

We consider a single river reach with 10 uniformly spaced water level nodes and rectangular cross section, an upstream inflow boundary condition provided with a fixed time series, as well as a controllable downstream release boundary condition. 
The grid is illustrated in Figure~\ref{fig:example-grid}, and the hydraulic parameters and initial conditions are summarized in Table~\ref{table:example-specs}. The model starts from steady state: the initial flow rate is uniform and the water level decreases linearly along the length of the channel.

\begin{table}[!t]
\captionsetup{font=scriptsize}
  \caption{Parameters for the example problem.}
  \label{table:example-specs}
  \centering
  \begin{tabularx}{\linewidth}{ l l X}
  \hline
     Parameter & Value & Description \\
  \hline
      $T$ & $72$ & Index of final time step \\
      $\Delta t$ & $600$ s & Time step size \\
      $H^b_{i}$ & $\left( -4.90, -4.92, \ldots, -5.10 \right)$ m & Bottom level\\
      $l$ & $10\,000$ m & Total channel length \\
      $A_i(H_i)$ & $50\cdot(H-H^b_i)$ m$^2$ & Channel cross section function \\
      $P_i(H_i)$ & $50+2\cdot(H-H^b_i)$ m & Channel wetted perimeter function \\
      $C_{i}$ & $\left( 40, 40, \ldots, 40 \right)$ m$^{0.5}$/s & Chézy friction coefficient \\
      $H_i(t_0)$ & $\left( 0.000, -0.025, \ldots, -0.222 \right)$ m & Initial water levels at $H$ nodes \\
      $Q_i(t_0)$ & $\left( 100, 100, \ldots, 100 \right)$ m$^3$/s & Initial discharge at $Q$ nodes \\
    \hline
  \end{tabularx}
  
\end{table}

Our objective is to keep the water level at the $H$ node upstream of the gate at $0$ m above datum:
\[
f=\sum_{j=1}^T |H_{10}(t_j)|^2,
\]
\emph{i.e.}, $f=\|y\|_2^2$ and $y=\left(H_{10}(t_j)\right)_{j=1,\ldots,T}$ in concordance with Conditions 7 and 8 of Definition \ref{def:rmpc},
subject to the bounds on the control variables
\begin{equation*}
100\, \text{m}^3/\text{s} \leq Q_{10}(t_j) \leq 200\, \text{m}^3/\text{s},
\end{equation*}
in concordance with Conditions 1--2 and 5 of Definition \ref{def:rmpc}.

\begin{figure}[htp]
  \centering

  \subfloat[Fixed inflow boundary condition.]{\label{fig:results1}\includegraphics[width=60mm]{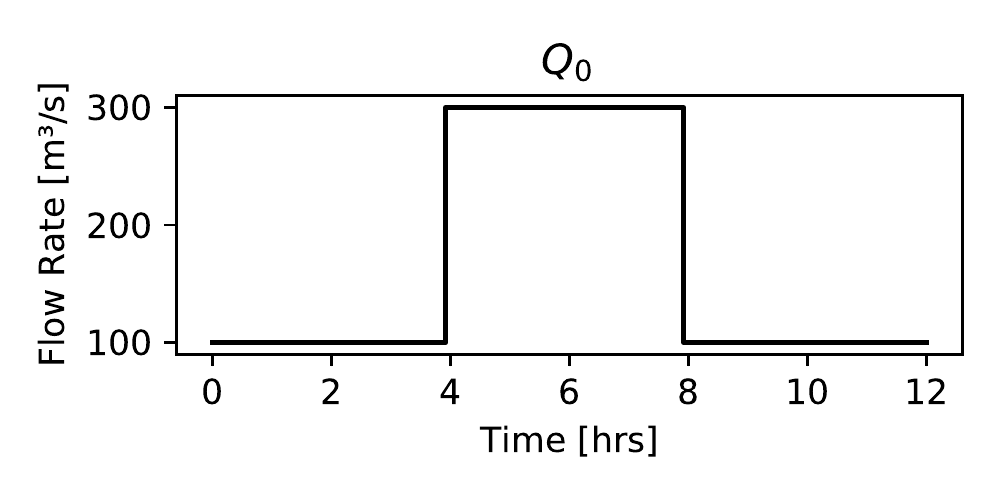}}
  \subfloat[Outflow bounds (dashed) and trajectory.]{\label{fig:results2}\includegraphics[width=60mm]{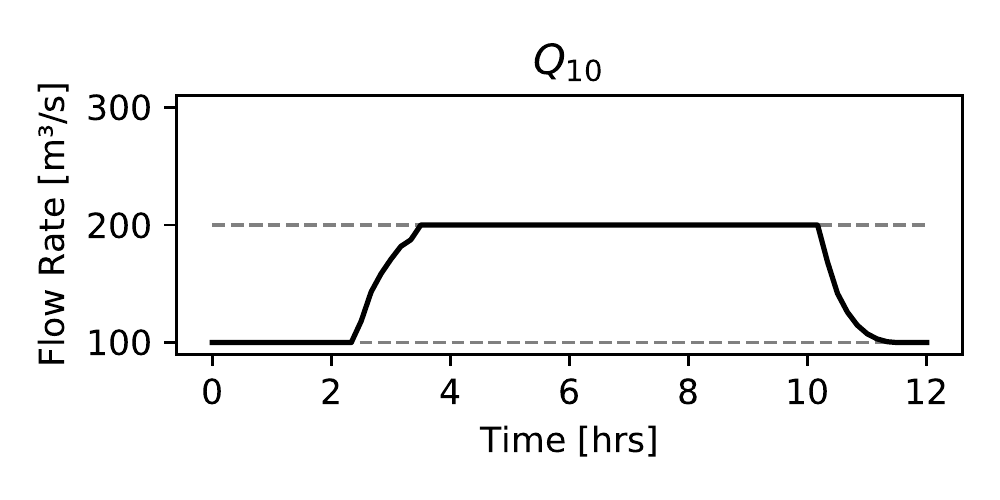}}
  \\
  \subfloat[Upstream level  trajectory.]{\label{fig:results3}\includegraphics[width=60mm]{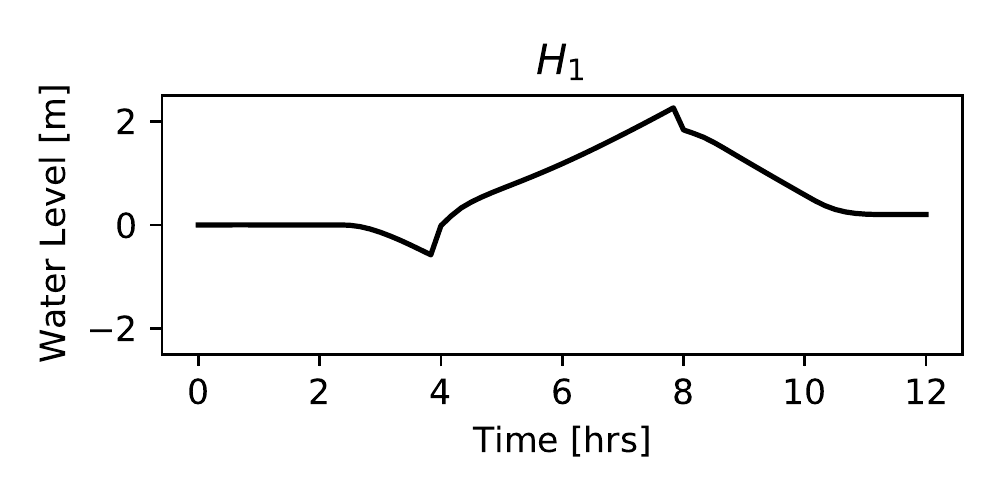}}
  \subfloat[Downstream level  trajectory.]{\label{fig:results4}\includegraphics[width=60mm]{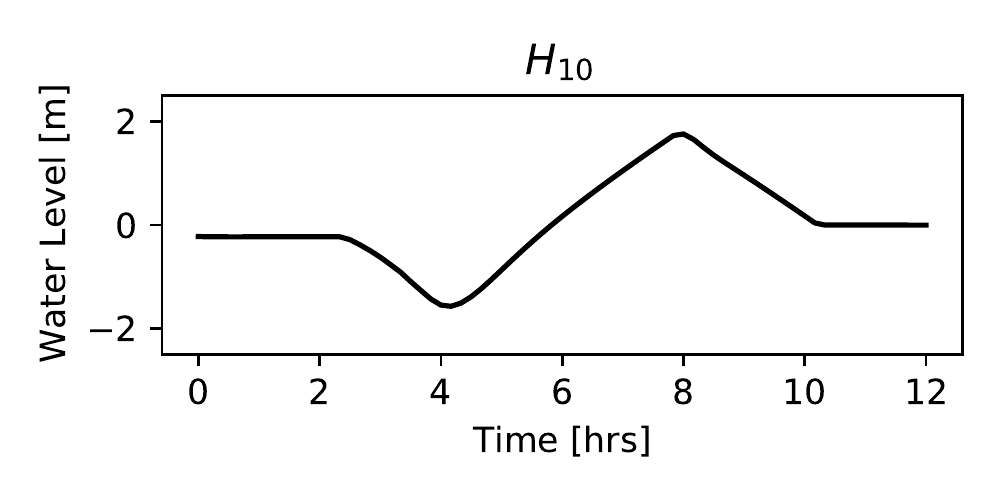}}
  \captionsetup{font=scriptsize}
  \caption{Solution to the example problem.}\label{fig:results}

\end{figure}

To give a physical context for this problem, suppose this model represents a channel downstream of a reservoir and upstream of an adjustable gate with limited capacity.
The gate is trying to dampen the sudden pulse of water shown in Figure~\ref{fig:results1} released by the reservoir.

A solution to the optimization problem was obtained using the interior point solver \textsc{IPOPT} \cite{Wachter2006} and is plotted in Figure~\ref{fig:results}. By releasing water in anticipation of the inflow using the decision variable $Q_{10}$, the optimization is able to reduce water level fluctuations and keep the water levels close to the target level.

The \textsc{IPOPT} solve, starting from an all-zero starting point, takes approximately $0.1$ s to complete on a $2.6$ GHz Intel Core i7 CPU.

Since some of the bounds on the control variables are active, the near-global optimality guarantee of Theorem \ref{thm:invex} applies. The (near-)global optimality was checked numerically by seeding the optimization with a large number of different starting points. Latin hypercube sampling \cite{mckay2000comparison} was used to compute $1\,000$ different starting points, for each of which \textsc{IPOPT} computed a solution. The standard deviation of the solution vectors was found to be in the order of $10^{-11} - 10^{-13} \approx 0$ per  solution vector coordinate, illustrating how every starting point resulted in an -- for all practical purposes -- identical solution.  This provides evidence  that the found solution is indeed globally optimal, or at least nearly so.

There is also ample other numerical evidence that solutions of this type are globally optimal or very close to it.  In \cite{Baayen2019-3}, the performance of an interior point-type method (IPM) for a large class of water problems is benchmarked against a so-called \emph{reduced genetic algorithm} (RGA, \cite{vermuyten2018combining}).
The IPM search finds qualitatively consistent solutions that always obtains better objective function values than the RGA. This benchmark includes problems with multiple river reaches, multiple spatial control points, and both coarser and finer discretizations of the shallow water equations in time and space.

Similar results are reported in \cite{ghaddar2017polynomial} for drinking water distribution networks, where local search using \textsc{IPOPT} finds solutions with objective values within a relative distance of $10^{-3}$ of those found using the global solver \textsc{Couenne} -- in a fraction of the computation time. The general-purpose global solvers require multiple minutes or hours to run, but the local search completes in a few seconds at most.

\section*{Acknowledgements}
The authors would like to thank Dick den Hertog and Vyacheslav Kungurtsev for the critical reading of the manuscript and their comments. The second author's work was financed by The Dutch Research Council (NWO) grant VI.Veni.191E.035.

\printbibliography

\end{document}